\documentclass{article}
\def\n{\noindent}
\begin{document}

\title {\bf On the divergence of greedy algorithms with respect to Walsh subsystems in $L$}
\author{S.A.Episkoposian}

\date {e-mail: sergoep@ysu.am}

\maketitle {In this paper we prove that there exists a function
which $f(x)$ belongs to $L^1[0,1]$ such that a greedy algorithm
 with regard to the Walsh subsystem does not converge to $f(x)$
in $L^1[0,1]$ norm, i.e. the Walsh subsystem $\{W_{n_k}\}$ is not
a quasi-greedy basis in its linear span in $L^1$ } \vskip 4mm

{\bf 1. INTRODUCTION }
\par\par\bigskip

Let a Banach space $X$ with a norm $||\cdot||=||\cdot||_X$, and a
basis $\Psi=\{\psi_k\}_{k=1}^\infty$, $||\psi_k||_X=1$, $k=1,2,..$
be given.

Denote by $\Sigma_m$ the collection of all functions in $X$ which
can be expressed as a linear combination of at most $m$- functions
of $\Psi$. Thus each function $g\in \Sigma_m$ can be written in
the form
$$g=\sum_{s\in \Lambda}a_s\psi_s,\ \ \#\Lambda \leq m.$$

For a function $f\in X$ we define its approximation error
$$\sigma_m(f,\Psi)=\inf_{g\in \Sigma_m}\left| \left| f-g \right| \right|_X,\ \ m=1,2,...$$

and we consider the expansion
$$f=\sum_{k=1}^\infty a_k(f)\psi_k.$$

{\bf Definition 1.} Let an element $f\in X$ be given. Then the
$m$-th greedy approximant of function $f$ with regard to the basis
$\Psi$ is given by the formula

$$G_m(f,\Psi)=\sum_{k\in\Lambda}a_k(f) \psi_k,$$

where $\Lambda \subset \{1,2,...\}$ is a set of cardinality $m$
such that $$|a_n(f)|\geq |a_k(f)|, \ \  n\in \Lambda,\ \  k\notin
\Lambda.$$
\vskip 4mm

We'll say that the greedy approximant of $f(t)\in L^p_{[0,1]}$,
$p\geq 1$ with regard to the basis $\Psi$ converges, if the
sequence $G_m(x,f)$ converges to $f(t)$ in $L^p$ norm.

This new and very important direction invaded many mathematician's
attention (see [1]-[9]).

{\bf Definition 2}. We call a basis $\Psi$ greedy basis if for
every $f \in X$  there exists a subset $\Lambda\subset
\{1,2,...\}$ of cardinality $m$, such that
$$\left| \left| f-G_m(f,\Psi) \right| \right|_X\leq C\cdot\sigma_m(f,\Psi)$$
where a constant $C=C(X,\Psi)$ independent of $f$ and $m$. \vskip
4mm

In 1998 V.N.Temlyakov proved that each basis $\Psi$ which is $L_p$
-equivalent to the Haar basis $H$ is Greedy basis for $L_p(0,1)$,
$1<p<\infty$ (see [4]).

\vskip 4mm

{\bf Definition 3.} We say that a basis $\Psi$ is Quasi-Greedy
basis if there exists a constant $C$  such that for every $f\in X$
and any finite set of indices $\Lambda$, having the property
$$\min_{k\in \Lambda}|a_k(f)|\geq \max_{n\notin \Lambda}|a_k(f)|$$
we have
$$\left| \left| \sum_{k\in \Lambda}a_k(f)\psi_k \right| \right|_X \leq C\cdot ||f||_X.$$
\vskip 4mm

In 2000 P.Wojtaszczyk [5] proved that a basis $\Psi$ is
quasi-greedy if and only if the sequence $\{G_m(f) \}$ converges
to $f$, for all $f\in X$. Note that in [6] S.Konyagin and
V.Temlyakov constructed an example of quasi-greedy basis that is
not Greedy basis.

V.Temlyakov proved that the trigonometric system $T$ is not a
Quasi-Greedy basis for $L^p$ if $p\not=2$ (see [7]).

In [8] it is proved that this result is true for Walsh system.

In the sequel, we'll fix a sequence $\{ M_n \}_{n=1}^\infty$ so
that
$$\lim_{k \to \infty}(M_{2k}-M_{2k-1})=+\infty$$
and consider a subsystem of Walsh system
$$\{ W_{n_k}(x) \}_{k=1}^\infty=\{W_ m(x):\ \ M_{2s-1}\leq m \leq M_{2s},\ \ s=1,2,...\}\eqno(1)$$

 In this paper we constructed a function $f(x)\in
L^1[0,1]$ such that the sequence $\{G_m(f) \}$, with respect to
Walsh system, does not converge to $f(x)$ by $L^1[0,1]$ norm and
we can watch for spectra of "bad" function $f(x)$.

Moreover the following is true.

 \vskip 4mm

{\bf Theorem .}  There exists a function $f(x)$ belongs to
$L^1[0,1]$ such that the approximate $G_n(f,W_{n_k})$ with regard
to the Walsh subsystem does not converge to $f(x)$ by $L^1[0,1]$
norm, i.e. the Walsh subsystem $\{W_{n_k}\}$ is not a quasi-greedy
basis in its linear span in $L^1$.

\par\par\bigskip
\par\par\bigskip

{\bf 2. PROOF OF THEOREM }
\par\par\bigskip
\par\par\bigskip

First we will give a definition of Walsh-Paly system (see [10]).
$$W_0(x)=1,\ \ {\displaystyle W_n(x)= \prod_{s=1}^k r_{m_s}(x),\ \ n= \sum_{s=1}^k 2^{m_s}},\ \ m_1>...>m_s, \eqno(2)$$
where $\{ r_k(x)\}_{k=0}^\infty$ is the system of Rademacher:
$$r_0(x)=\cases{\hskip8pt 1,\ \ x \in \left[0,{1\over2}\right); \cr \\ \cr -1, \ \ x \in ({1\over2},1).\cr} $$
$$r_0(x+1)=r_0(x),\ \ r_k(x)=r_0(2^kx),\ \ k=1,2,...$$

In the proof of theorem we will used the following properties of
Walsh system:

\vskip 4mm

1. From (2) we have
$$W_{2^k+j}(x)=W_{2^k}(x) \cdot W_j(x),\ \ if\ \ 0 \leq j \leq 2^k-1.\eqno (3)$$

2. The Dirichlet-Walsh kernel $\displaystyle
D_m(x)=\sum_{j=0}^{m-1}W_j(x)$ has the following properties (see
[10] p.27)

$$D_{2^j}(x)=\cases{2^j , \ \ x \in \big[0,{1\over 2^j}\big];\cr \cr 0 ,\ \ x \in \big({1\over 2^j},1\big].\cr}\eqno (4)$$

3. There is a sequence of natural numbers $\{m_k\}_{k=1}^\infty$
with $2^{k-1}\leq m_k <2^k$, $k=1,2,...$ (see [10] p.47), such
that
$$||D_{m_k}||_1=\int_0^1|D_{m_k}(x)|dx \geq {1\over 4}\log_2 m_k,\ \ k=1,2,.... \eqno(5)$$

{\bf Proof of Theorem .}  Taking into account (1)-(3) we can take
the sequences of natural numbers $\{k_\nu \}_{\nu=1}^\infty$ and
$\{p_\nu \}_{\nu=1}^\infty$ so that the following conditions are
satisfied:

$$k_\nu>(\nu-1)^2+1, \eqno(6)$$
$$W_{2^{k_\nu}}(x) \cdot W_i(x)=W_{2^{k_\nu}+i}(x),\ \  0 \leq i < 2^{k_\nu},\eqno (7)$$
$$2^{k_\nu}+i\in [ M_{2p_\nu-1},M_{2p_\nu}),\ \  0 \leq i < 2^{k_\nu},\eqno (8)$$

For any natural $\nu$ we set
$$f_\nu(x)=\sum_{N_{\nu}\leq n_k <N_\nu-1} c_{n_k}^{(\nu)}W_{n_k}(x)=$$
$$=\sum_{i=0}^{2^{k_\nu}-1}\left( {1 \over \nu^2}+2^{-(2^{k_\nu}+i)}\right) \cdot W_{2^{k_\nu}+i}(x)=$$
$$=W_{2^{k_\nu}}(x)\cdot \sum_{i=0}^{2^{k_\nu}-1}\left( {1 \over \nu^2}+2^{-(2^{k_\nu}+i)}\right) \cdot W_i(x)=$$
$$=W_{2^{k_\nu}}(x)\cdot \left[ {1 \over \nu^2}\sum_{i=0}^{2^{k_\nu}-1}W_i(x)+ {1 \over 2^{k_\nu}}\sum_{i=0}^{2^{k_\nu}-1}2^{-i}W_i(x)
\right]=$$
$$=W_{2^{k_\nu}}(x)\cdot \left[ {1 \over \nu^2}D_{2^{k_\nu}}(x)+ {1 \over 2^{k_\nu}}\sum_{i=0}^{2^{k_\nu}-1}2^{-i}W_i(x)\right],\eqno(9)$$
where
$$c_{n_k}^{(\nu)}=\cases{{1 \over \nu^2}+2^{-n_k}, \ \ N_{\nu} \leq n_k =N_{\nu}+i< N_{\nu+1}, \ \ 0 \leq i<N_{\nu}, \cr \\  \cr   0 ,\ \ n_k<N_{\nu},\ \ \nu\geq 1, \cr}\ \  \eqno(10)$$
$$N_{\nu}=2^{k_\nu},\ \ N_{\nu+1}=2^{k_\nu+1}. \eqno(11)$$

We set
$$f(x)=\sum_{k=1}^\infty c_{n_k}(f)W_{n_k}(x)=$$
$$=\sum_{\nu=1}^\infty f_\nu(x)=\sum_{\nu=1}^\infty \left[ \sum_{N_{\nu}\leq n_k <N_\nu-1} c_{n_k}^{(\nu)}W_{n_k}(x)\right], \eqno(12) $$
where
$$c_{n_k}(f)=c_{n_k}^{(\nu)} \ \ for \ \ N_{\nu}\leq n_k <N_\nu-1,\ \ \nu=1,2,...\eqno(13)$$
Now we will show that $f(x)\in L^1[0,1]$

Taking into account (9)-(11)we get
$$f(x)=\sum_{\nu=1}^\infty {1 \over \nu^2}W_{2^{k_\nu}}(x)D_{2^{k_\nu}}(x)+\sum_{\nu=1}^\infty \left[\sum_{i=0}^{2^{k_\nu}-1}2^{-i}W_{2^{k_\nu}+i}(x) \right]=$$
$$=G(x)+H(x).\eqno(14)$$

For function $G(x)$ from (4)and definition of Walsh system we have
$$\int_0^1|G(x)|dx\leq \sum_{\nu=1}^\infty {1 \over \nu^2}<\infty$$
and we get that $G(x)\in L^1[0,1]$.

Analogously
$$\int_0^1|H(x)|dx\leq \sum_{\nu=1}^\infty 2^{-\nu}<\infty$$
i.e. $H(x)\in L^1[0,1]$.

Hence and from (14) it follows that $f(x)\in L^1[0,1]$.

For any natural $\nu$ we choose numbers $k,j$ so that
$$N_{\nu}\leq n_k <N_{\nu+1}\leq n_j <N_{\nu+2}.$$

Then from (10) we have
$$c_{n_j}(f)=c_{n_j}^{(\nu+1)}={1 \over (\nu+1)^2}+2^{-n_j}<$$
$$<{1 \over \nu^2}+2^{-n_k}=c_{n_k}^{(\nu)}=c_{n_k}(f)$$
i.e. $c_{n_j}(f)<c_{n_k}(f)$.

Analogously for any number $n_k$, $N_{\nu}\leq n_k <N_{\nu+1}$,
$\nu \geq 1$ we have
$$c_{n_k+1}^{(\nu)}={1 \over \nu^2}+2^{(-n_k+1)}<{1 \over \nu^2}+2^{-n_k}=c_{n_k}^{(\nu)}$$

Thus we get
$$c_{n_{k+1}}(f)<c_{n_k}(f).$$

In other hand if $k\to \infty$ then $n_k\to \infty$ and $\nu \to
\infty$ (see (10), (11)).

From (13) we get $\displaystyle \lim_{k\to \infty}c_{n_k}(f)=0$
and consequently $c_{n_k}(f) \searrow 0$.

For any numbers $m_\nu $ so that
$$2^{k_\nu} \leq m_\nu<2^{k_\nu+1}.\eqno(15)$$

By (11) - (13) and Definition 1 we have

$$G_{2^{k_\nu}+m_\nu}(f,W_{n_k})-G_{2^{k_\nu}}(f,W_{n_k})= $$
$$=\sum_{i=2^{k_\nu}}^{2^{k_\nu}+m_\nu-1}\left( {1 \over \nu^2}+2^{-(2^{k_\nu}+i)}\right) \cdot W_{2^{k_\nu}+i}(x)= $$
$$={1 \over \nu^2}\cdot W_{2^{k_\nu}}(x)\cdot \sum_{i=2^{k_\nu}}^{2^{k_\nu}+m_\nu-1}W_i(x)+$$
$$+ {1 \over 2^{k_\nu}}\cdot W_{2^{k_\nu}}(x)\cdot\sum_{i=2^{k_\nu}}^{2^{k_\nu}+m_\nu-1} {1 \over 2^i}W_i(x)=$$
$$=J_1+J_2. \eqno(16)$$

By (6) we get
$$J_1={1 \over \nu^2}\cdot W_{2^{k_\nu}}(x)\cdot \sum_{i=0}^{m_\nu-1}W_{2^{k_\nu}+i}(x)=$$
$$={1 \over \nu^2}\cdot W^2_{2^{k_\nu}}(x)\cdot D_{m_\nu}(x). $$

$$\mid J_2 \mid \leq \sum_{i=2^{k_\nu}}^{2^{k_\nu}+m_\nu-1} {1 \over 2^i} \mid W_i(x)\mid \leq \sum_{i=2^{k_\nu}}^{\infty} {1 \over 2^i}\leq 2^{2^{-k_\nu}+1}.$$

From this and (16)we obtain
$$\mid G_{2^{k_\nu}+m_\nu}(f,W_{n_k})-G_{2^{k_\nu}}(f,W_{n_k}) \mid \geq $$
$$\geq {1 \over \nu^2}\cdot \mid D_{m_\nu}(x)\mid-2^{2^{-k_\nu}+1}.\eqno(17)$$

Now we take the sequence of natural numbers $m_\nu$ defined as (5)
such that $2^{k_\nu}\leq m_\nu < 2^{k_\nu+1}$.
 Then from (6), (17) we have
$$\int_0^1\mid G_{2^{k_\nu}+m_\nu}(f,W_{n_k})-G_{2^{k_\nu}}(f,W_{n_k})\mid dx>$$
$$>{1 \over \nu^2}\cdot \int_0^1|D_{m_\nu}(x)|dx - 2^{2^{-k_\nu}+1} \geq $$
$$\geq{1\over 4\cdot \nu^2}\cdot \log_2 m_\nu - 2^{2^{-k_\nu}+1} \geq  {k_\nu \over 4 \cdot \nu^2} - 2^{2^{-k_\nu}+1} \geq $$
$$\geq {(\nu-1)^2+1 \over 4 \cdot \nu^2} - 2^{2^{-k_\nu}+1}\geq {1 \over 8}- 2^{2^{-k_\nu}+1} \geq C_1,\ \ \ \ \nu \geq 2$$

Thus  the sequence $\{G_n(f,W)\}$ does not converge  by $L^1[0,1]$
norm, i.e. the Walsh subsystem $\{W_{n_k}\}_{k=1}^\infty$ is not a
quasi-greedy basis in its linear span in $L^1$.

\par\par\bigskip
{\bf The Theorem is proved.}
\par\par\bigskip

{\bf Remark.} As we well known (see [10] p.149) if the
$c_i\searrow 0$, then the series $\displaystyle \sum_{n=1}^\infty
c_nW_n(x)$ converges on $(0,1)$. In the proof of Theorem we
constructed the series (11) so that the coefficients strongly
decreasing, but the series diverges by $L^1$ -norm.

\par\par\bigskip
\par\par\bigskip
\break

 {\bf R E F E R E N C E S }

\par\par\bigskip
\par\par\bigskip
\par\par\bigskip

\n [1] DeVore R. A., Temlyakov V. N., Some remarks on Greedy
Algorithms, Adv. Comput. Math., 1995, v.5, p.173-187.
\par\par\bigskip

\n [2] DeVore  R. A., Some remarks on greedy algorithms. Adv
Comput. Math. 5, 1996, 173-187.
\par\par\bigskip

\n [3] Davis G., Mallat S. and Avalaneda M. 1997, Adaptive greedy
approximations. Constr, Approx. 13, 57-98.
\par\par\bigskip

\n [4] Temlyakov V. N., The best $m$ - term approximation and
Greedy Algorithms,  Advances in Comput. Math., 1998, v.8,
p.249-265.
\par\par\bigskip

\n [5] Wojtaszcyk P., Greedy Algorithm for General Biorthogonal
Systems, Journal of Approximation Theory, 2000, v.107, p. 293-314.
\par\par\bigskip

\n [6] Konyagin S. V., Temlyakov V. N., A remark on Greedy
Approximation in Banach spaces, East Journal on Approximation,
1999, v.5, p. 493-499.
\par\par\bigskip

\n [7] Temlyakov V. N., Greedy Algorithm and $m$ - term
Trigonometric approximation, Constructive Approx., 1998, v.14, p.
569-587.
\par\par\bigskip

\n [8] Gribonval R., Nielsen M., On the quasi-greedy property and
uniformly bounded orthonormal systems ,
http://www.math.auc.dk/research/reports/R-2003-09.pdf.
\par\par\bigskip

\n [9] Grigorian M.G., "On the convergence of Greedy algorithm",
International Conference, Mathematics in Armenia, Advances and
Perspectives, Abstract, 2003,  p.44 - 45, Yerevan, Armenia.
\par\par\bigskip

\n [10] Golubov B. I., Efimov A. V., Skvortsov V. A., Walsh Series
and Transformations: Theory and Applications [in Russian], NAuka,
Moscow, (1987); English transl.: Kluwer, Dordrecht (1991).
\par\par\bigskip
\par\par\bigskip
\par\par\bigskip

Department of Physics,

State University of Yerevan,

Alex Manukian 1, 375049 Yerevan, Armenia

e-mail: sergoep@ysu.am

\end{document}